
 \input amstex
\documentstyle{amsppt}
\TagsOnRight
\NoRunningHeads

\topmatter
\title
Invariant Manifolds of Hypercyclic Vectors for the Real scalar case.
\endtitle

\author
Juan P. B\`es
\endauthor
\address
DEPARTMENT OF MATHEMATICS AND COMPUTER SCIENCE, KENT STATE UNIVERSITY.
KENT, OH 44242
\email
jbes\@mcs.kent.edu
\endemail
\endaddress
\abstract
We show that every hypercyclic operator on a real locally convex
vector space admits a dense invariant linear manifold of
hypercyclic vectors.
\endabstract

\email jbes@mcs.kent.edu \endemail
\subjclass Primary 47A15 47A99 \endsubjclass
\thanks
The author wishes to thank the support
of the
Center for International and Comparative Programms 
and the Graduate Student Senate of
Kent State University.
\endthanks

\endtopmatter

\document

\define \w                      {\overline{w} }

\define \matrixad               {\bmatrix a_2 \\ b_2 \endbmatrix }

\define \matrixan               {\bmatrix a_n \\ b_n \endbmatrix }

\define \matrixo               {\bmatrix \format \r & \r \\ w & 0  \\ 0 & -w \endbmatrix }

\define \matrixu               {\bmatrix \format \r & \r \\ 1 & 1 \\ -\w & -w \endbmatrix }

\define \matrixw              {\bmatrix  \format \r & \r \\ -w & -1 \\ \w & 1 \endbmatrix }

Given a locally convex vector space $X$ and a continuous operator
$T:X\longrightarrow X$, we say that T is {\bf hypercyclic } provided there
exists some $x\in X$ whose orbit
$$
Orb\{ x , T \} = \{ x, Tx , T^2x, \dots \}
$$
is dense in $X$.\, Such an $x$ is said to be a hypercyclic vector 
for T.\,A motivation for this definition comes from the invariant subset
problem:
T has no non-trivial closed invariant subset if and only if every
non-zero vector in $X$ is hypercyclic for T.
\medskip
In 1990, B. Beauzamy constructed an example of an operator on a separable
complex Hilbert space admitting a dense, invariant linear manifold whose
non-zero vectors were all hypercyclic \cite{3, Thm.\, A}. Using different techniques, 
G. Godefroy and J. Shapiro provided a large class of examples on Fr\`echet 
spaces of entire functions with the same property \cite{5, Thm 5.1}. 
\smallskip
Soon after, D. Herrero \cite{6, Proposition 4.1} 
and P. Bourdon \cite{4} independently showed that {\it every} \  hypercyclic 
operator on a complex Hilbert space admits such a dense, invariant 
linear manifold of hypercyclic vectors (in fact, Bourdon's proof works for 
arbitrary complex locally convex spaces as well).  
We'd like to show here that the same 
holds for the real scalar case, by presenting a positive answer to the
following question, raised by S. Ansari \cite{1, Problem 1}:

\smallskip
``Suppose $X$ is a locally convex real vector space and $T:X\rightarrow X$ is
a continuous linear operator with a hypercyclic vector $x$. Is it true that
$P(T)x$ is a hypercyclic vector for $T$ whenever P is a non-zero polynomial
with real coefficients?''
\smallskip
(If so, 
$$\Cal M = \left \{ P(T)x \, : P \ \text{polynomial with real coefficients} \  
\right \} 
$$
 is a dense, $T$-invariant manifold of hypercyclic vectors).  
\smallskip
Notice that since $P(T)$ and $T$ commute,
$$
Orb\{ P(T)x , T \} = P(T) \left ( Orb\{ x , T \} \right ).
$$
That is, given x a hypercyclic vector for $T$, $P(T)x$ will also be
hypercyclic if and only if $P(T)$ has dense range. So it will suffice for 
us to show the following.

\proclaim{Theorem} Let $X$ be a locally convex real vector space, 
and $T\in L(X)$ be hypercyclic.
If $P$ is a non-zero polynomial with real coefficients, 
then $P(T)$ has dense range.

\endproclaim

We'll make use of two results. The first one is due to C.Kitai
\cite{6, Theorem 2.3}:

\proclaim{Lemma 1}
Let $X$ be a locally convex real ( complex ) vector space, and $T\in L(X)$ 
be  hypercyclic.
Then the adjoint $T^*$ of $T$ has no eigenvalues. In other words, for any
scalar $b$, the operator $T+bI$ must have dense range.
\endproclaim

\demo {Proof} \enddemo

Suppose $T^*$ admits an eigenvector $x$ with eigenvalue $\lambda$,
and let $z\in X$ be a hypercyclic vector for $T$. Since $x\ne 0$,
the set
$$
\left \{  \langle x , T^n z \rangle \right \}_{n\ge 1} =
\left \{  \langle {T^*}^n x , z \rangle \right \}_{n\ge 1}=
\left \{ {\lambda}^n  \langle x , z \rangle \right \}_{n\ge 1}
$$
must be dense in the real (complex) scalar field, contradiction. $\square $

\smallskip
\smallskip
The statement of the 
second result may be traced back at least to S. Rolewicz \cite{8, p17}.
We are grateful to P. Bourdon for suggesting the following argument to us.
\smallskip
\proclaim{Lemma 2}
$R^n$ admits no hypercyclic operators.
\endproclaim

\demo {Proof} \enddemo

Suppose there exists $A:R^n \rightarrow R^n$  linear with hypercyclic vector $z$.
Notice that $B=\{ z, Az, \dots , A^{n-1}z \}$ must be linearly independent.
Now, since $z$ is a hypercyclic vector for $A$, there exist sequences of 
positive integers $(n_k)$ and $({\tilde n}_k )$ so that
$$
\align
A^{n_k}z & \underset{k\to \infty}\to \rightarrow 0 \\ 
A^{\tilde{n}_k}z &\underset{k\to \infty}\to \rightarrow z. 
\endalign
$$
Because $B$ is a basis of $R^n$, we have in fact that
$$
\align
A^{n_k}x & \underset{k\to \infty}\to \rightarrow 0 \\                
A^{\tilde{n}_k}x &\underset{k\to \infty}\to \rightarrow x
\endalign  
$$
for all $x\in R^n$. Hence, if $|A|$ denotes the determinant of $A$, the
above lines imply the contradictory fact that
$$
\align
|A|^{n_k}&=|A^{n_k}|\underset{k\to \infty}\to \rightarrow 0 \\
|A|^{\tilde{n}_k}&=|A^{\tilde{n}_k}|\underset{k\to \infty}\to \rightarrow 1. \ \ 
\endalign 
$$
So Lemma 2 holds. \   \text{ $\square $}
\medskip

Now, let's show the Theorem.

\newpage 

\demo {Proof of Theorem } \enddemo

\smallskip
Since scalar multiples and compositions of operators having dense range have
dense range, we may assume P is irreducible and monic. Moreover, by Lemma 1
we may assume $P$ is of the form
$$
P(t)=t^2-2 Re(w)t +|w|^2 \ , \ \ \ \text{for some non-real complex number $w$.}
$$

Now, suppose that $P(T)$ does {\it not} have dense range. Let
$0\ne x \in Ker(P(T)^*)$. Then
$$
{T^*}^2 x = a_2 T^* x + b_2 x \ , \text{where} \ \ 
 \left \{ \matrix
 a_2 &= &w + \w \\
 b_2& = & -|w|^2 
 \endmatrix \right. . \tag1
$$

By Lemma 1, $\left \{ T^* x , x \right \} \subset Ker(P(T)^*)$ must be linearly 
independent.
So  there exist unique scalars $a_n$ and $b_n$ satisfying
$$
{T^*}^n x = a_n  T^* x + b_n x \ \ \ \ \ \ \ \ (\, n=1,2,\dots\, ).  
\tag 2
$$

Notice that by (2)
$$
\align
{T^*}^{n+1}x &= T^* (a_n T^*x +b_n x) \\
             &= a_n (a_2 T^*x +b_2 x) + b_n T^*x \\
             &= (a_n a_2 + b_n ) T^*x + a_n b_2 x .
\endalign
$$

That is, 
$$
\bmatrix a_{n+1} \\
         b_{n+1}
\endbmatrix
=
A \bmatrix a_{n} \\
b_n
\endbmatrix
\ \  \ , \text{where} \ \ \  A =
\bmatrix
a_2 & 1 \\
b_2 & 0
\endbmatrix
=
\bmatrix
w+\overline{w} & 1 \\
-|w|^2 & 0
\endbmatrix .
$$
\smallskip
So for all $n\ge 2$,

$$
\matrixan
= A^{n-2}
\matrixad
 \tag 3
$$
and 

$$
\bmatrix
a_2 & b_2 \\
a_3 & b_3
\endbmatrix
\ 
\left ( A^{n-2} \right )^{t} 
=
\bmatrix
a_n & b_n \\
a_{n+1} & b_{n+1}
\endbmatrix .
\tag 4
$$
\medskip

Next, notice that since 
$\left \{ T^* x , x \right \}$ is linearly independent, the operator
$$
\align
&{\CD
X @ > \Pi >> R^2  
\endCD } \\
& y \mapsto ( \langle x , y \rangle ; \langle T^* x ,y\rangle \, )
\endalign
$$
is onto. In particular, if $z\in X$ is a hypercyclic vector for $T$, then
$$
\left \{ \ \Pi ( T^n z) \ \right \}_{n\ge1} \tag 5
$$
must be dense in $R^2$.

\smallskip
Hence, by (2) and (4) we have for each $n\ge 2$
$$
\spreadlines{1\jot }
\align
\Pi ( T^n z)  &= ( \langle x , T^n z \rangle ; \langle T^* x ,T^n z \rangle ) \\
	      &= (  \langle {T^*}^n x , z \rangle ; \langle {T^*}^{n+1} x , z \rangle ) \\
	      &= (  \langle a_n T^* x + b_n x , z \rangle ; \langle a_{n+1} T^* x + b_{n+1}
	       x , z \rangle ) \\
	      &= {\bmatrix a_n & b_n \\ a_{n+1} & b_{n+1} \endbmatrix}
	          {\bmatrix \langle T^*  x, z\rangle \\ \langle x , z \rangle \endbmatrix}
		   \\
	      &= {\bmatrix a_2 & b_2 \\ a_3 & b_3 \endbmatrix} \left ( A^{n-2}\right )^t 
	         {\bmatrix \langle T^*  x, z\rangle \\ \langle x , z \rangle \endbmatrix} .
\tag6
\endalign
$$
\medskip
Thus, (5) and (6) force $A^t$ to be hypercyclic on $R^2$, which contradicts Lemma 2.

So $P(T)$ must have dense range.\ \ $\square$

\remark {Remark}

As a consequence of the theorem, S. Ansari's proof that every operator
on a complex Banach Space shares with its powers the {\it same }
hypercyclic vectors \cite{2, Thm 1} works for the real scalar case as well
(See also \cite{1, Note 3}).
 
\endremark

\enddocument